\documentclass[smallextended]{svjour3}       
\smartqed  
\usepackage{graphicx}
\usepackage{mathptmx}      
\usepackage{amsmath}
\usepackage{amsfonts}
\usepackage{amssymb}

\usepackage{bm}
\DeclareMathAlphabet{\mathcal}{OMS}{cmsy}{m}{n}
\usepackage{hyperref}
\usepackage{mathtools}
\mathtoolsset{showonlyrefs=true}

\usepackage{algorithm}
\usepackage{algorithmic}

\usepackage{color}
\usepackage{dcolumn}

\usepackage{listings}

\newcommand{\trans}{{\text{\tiny\sf T}}}

\newcommand{\bbR}{\mathbb{R}}

\newcommand{\rmd}{\mathrm{d}}

\newcommand{\veczero}{\bm 0}
\newcommand{\vecone}{\bm 1}

\newcommand{\vecb}{\bm b}
\newcommand{\vecc}{\bm c}
\newcommand{\vecd}{\bm d}

\newcommand{\vecx}{\bm x}
\newcommand{\vecy}{{\bm y}}

\newcommand{\dg}{\nabla_{\mathrm d}}

\DeclareMathOperator{\diag}{diag}

\DeclareMathOperator{\rank}{rank}

\usepackage{tikz}
\usepackage{pgfplotstable}
\usetikzlibrary{arrows,positioning,plotmarks,external,patterns,angles,
decorations.pathmorphing,backgrounds,fit,shapes,graphs,calc}
\begin{document}

\title{Adaptive projected SOR algorithms for nonnegative quadratic programming}
\subtitle{}


\author{Yuto Miyatake \and Tomohiro Sogabe
}


\institute{
  Y. Miyatake \at
  Cybermedia Center, Osaka University, Osaka, Japan\\
  \email{miyatake@cas.cmc.osaka-u.ac.jp}
  \and 
  T. Sogabe \at
  Department of Applied Physics, Graduate School of Engineering, 
  Nagoya University, Nagoya, Japan
  \email{sogabe@na.nuap.nagoya-u.ac.jp}
}

\date{Received: date / Accepted: date}

\maketitle

\begin{abstract}
  The choice of relaxation parameter in
  the projected successive overrelaxation (PSOR) method for nonnegative quadratic programming problems is problem-dependent. 
  We present novel adaptive PSOR algorithms that adaptively control the relaxation parameter using the Wolfe conditions.
  The method and its variants can be applied to various problems without requiring additional assumptions, 
  barring the positive semidefiniteness concerning the matrix that defines the objective function,
  and the cost for updating the parameter is negligible in the whole iteration.
  Numerical experiments show that the proposed methods often perform comparably to (or sometimes superior to) the PSOR method with a nearly optimal relaxation parameter.

\end{abstract}

\section{Introduction}
\label{sec1}

Given a symmetric positive semidefinite matrix $A\in\bbR^{n\times n}$ and vector $\vecb\in\bbR^n$,
we consider the following convex quadratic programming problem:
\begin{equation}
  \tag{P1} \label{nqp}
  \begin{split}
    \text{minimize} & \quad V(\vecx) = \frac{1}{2} \vecx^\trans A \vecx - \vecx^\trans \vecb \\
    \text{subject to} & \quad  \vecx \geq \veczero \quad (x_i\geq 0 \text{ for all } i).
  \end{split}
\end{equation}
Throughout this paper, we assume that $A$ has positive diagonal entries.
The problem~\eqref{nqp} is referred to as, for example, the nonnegative quadratic programming (NQP) or symmetric monotone linear complementarity problem.
The problem~\eqref{nqp} is related to the nonnegative constrained least squares problem~\cite{bj96}
\begin{equation}
  \tag{P2} \label{nls}
  \begin{split}
    \text{minimize} &\quad  V(\vecx) = \| C\vecx - \vecd\|_2^2 \\
    \text{subject to} & \quad  \vecx \geq \veczero,
  \end{split}
\end{equation}
where $C\in\bbR^{m\times n}$, $\vecd\in\bbR^{m}$, 
such that \eqref{nqp} is converted to \eqref{nls}
with $A = C ^ \trans C$ and $\vecb = C^\trans \vecd$
, and vice versa.

The problem~\eqref{nqp} or its equivalent formulation arises from numerous contexts, some examples of which are listed here. 
Early examples include, for example, a linear complementarity problem for solving free boundary problems for journal bearings by finite differences~\cite{ch41}. 
In the context of machine learning, the problem is useful for computing the maximum margin hyperplanes in support vector machines~\cite{bu98}.
Some image reconstruction or deblurring techniques also require solving this problem~\eqref{nls}~\cite{ca04,ns00,st08}.
The problem also appears in the study of financial derivatives~\cite{wh95} and
 dynamic multi-rigid-body contact simulations~\cite{ap97}.

In this study, we are concerned with solving \eqref{nqp} or \eqref{nls} numerically.
Iterative algorithms for solving \eqref{nqp} or \eqref{nls} have been studied for more than half a century.
A straightforward approach is to apply a gradient method
to $V(\vecx)$ with an appropriate projection.
Furthermore, sophisticated variants exist based on matrix splitting.
Some can be interpreted as a stationary iterative method such as the successive overrelaxation (SOR) method for $A\vecx = \vecb$ with an appropriate projection in each iteration~\cite{cr71a}.
Other stationary iterative methods such as the accelerated overrelaxation method with an appropriate projection have recently been studied~\cite{ht15,ht16}.
Multiplicative iteration methods~\cite{sh02,sh07}, which are highly parallelizable, and interior-point methods~\cite{co09} are also popular in the context of machine learning and linear complementarity problems, respectively.
A class of inner-outer iterative methods has also been studied~\cite{ol80,po69,zh16}.

This paper focuses on the SOR method with an appropriate projection, which is called the projected SOR (PSOR) method.
The SOR method can be understood as a cyclic coordinate descent method as shown in section~\ref{subsec:geom_sor}, and the PSOR method projects the iterate for each componentwise update.
This approach is expected to be effective for large scale problems witout requiring the estimation of some properties of the problem such as the largest absolute eigenvalue of $A$.
As is the case with the SOR method applied to a linear system $A\vecx = \vecb$, the convergence performance depends on the relaxation parameter $\omega$.
Notably, for \eqref{nqp}, as we will see through numerical experiments, the convergence behaviour often drastically changes during the iteration.
Therefore, changing the relaxation parameter during the iteration may be of great use.
We consider controlling the relaxation parameter adaptively and call such a method an adaptive projected SOR (APSOR) method.
We aim to develop APSOR methods with the following properties: 
they perform effectively for a wide variety of~\eqref{nqp} without requiring a specific property for $A$, and 
the cost for updating the parameter is negligible in the overall iterations.
We note that existing adaptive SOR techniques developed for solving linear systems can be incorporated with the PSOR method for \eqref{nqp}; however, most of them require a specific structure for matrix $A$ or take a substantial cost for updating the relaxation parameter (a brief review will be given in section~\ref{subsec:apsor}).
Recently, the present authors proposed an adaptive SOR (ASOR) method for $A\vecx=\vecb$ with the above two properties~\cite{ms20}.
In this paper, we show that the idea in~\cite{ms20} is also useful for~\eqref{nqp}.
The main contributions of the present study are as follows.

\begin{itemize}
  \item 
        We point out that the PSOR method has a form of projected gradient methods.
        This is a straightforward but essential consequence of our previous report~\cite{ms18}, where we showed that the SOR method can be viewed as a discrete gradient scheme applied to a certain gradient system.
        The key of the discussion is that the relaxation parameter can be viewed as the step size after an appropriate change of variables.      
  \item 
        Using the above interpretation, we present an APSOR method based on the Wolfe conditions, which can be viewed as an application of~\cite{ms20} to \eqref{nqp}.
        Preliminary numerical experiments suggest that the convergence behaviour may fluctuate because of the adaptive control of the relaxation parameter; the convergence behaviour strongly depends on the initial guess, which is particularly typical for an ill-conditioned $A$ or symmetric positive semidefinite $A$.
        To address these issues, we upgrade the APSOR method through two ways.
        First, we develop an algorithm that fixes the relaxation parameter after performing some iterations.
        Second, we propose an algorithm that finds a reasonable initial guess by shifting the matrix $A$ in \eqref{nqp} and applying the APSOR method to the shifted problem.
\end{itemize}

The paper is organized as follows.
In section~\ref{sec:psor}, the PSOR method is briefly reviewed.
New APSOR algorithms are presented in section~\ref{sec:apsor} and numerically tested in section~\ref{sec:numer}.
Finally, we provide concluding remarks in section~\ref{sec:conclude}.

\section{Projected SOR method}
\label{sec:psor}

The PSOR method for solving~\eqref{nqp} was developed more than a half-century ago~\cite{cr71a}.
In this section, we review some basics of the method.

The SOR method is a typical stationary iterative method solving a linear system $A\vecx = \vecb$~\cite{go13,va00}.
Its iteration formula $\vecx^{(k)} \mapsto \vecx^{(k+1)}$ is defined as follows:
\begin{equation}
  \label{sor}
  x_i^{(k+1)} = (1-\omega) x_k^{(k)} + \frac{\omega}{a_{ii}} \Big( b_i - \sum_{j<i} a_{ij} x_j^{(k+1)} - \sum_{j>i} a_{ij} x_j^{(k)}\Big), \quad i=1,\dots,m,
\end{equation}
where $\omega$ is the relaxation parameter.
When $A$ is symmetric positive definite, the iterates will converge to $A^{-1}\vecb$ for any $\omega$ in the open interval $(0,2)$ and for any initial guess $\vecx^{(0)}$.

The PSOR method for \eqref{nqp} was developed around 1970~\cite{cr71a}. 
Its iteration formula is defined as follows:
\begin{equation}
  x_i^{(k+1)} = \max \Big[ (1-\omega) x_k^{(k)} + \frac{\omega}{a_{ii}} \Big( b_i - \sum_{j<i} a_{ij} x_j^{(k+1)} - \sum_{j>i} a_{ij} x_j^{(k)}\Big),0\Big], \quad i=1,\dots,m.
\end{equation}
The algorithm is summarized in Algorithm~\ref{algo:sor}.
Note that in Algorithm~\ref{algo:sor}, convergence is checked by monitoring $\| \vecx^{(k+1)} - \vecx^{(k)} \|$ as a simple example of the convergence criterion.
Other criteria are also possible.

\begin{algorithm}[t]
  \caption{PSOR method for~\eqref{nqp}.}
  \begin{algorithmic}[1]
    \label{algo:sor}
    \REQUIRE{$\omega \in (0,2)$, $\vecx^{(0)} \geq \veczero$}
    \FOR{$k=0,1,2,\dots$ \textbf{until} $\| \vecx^{(k+1)} - \vecx^{(k)} \| \leq \epsilon $}
    \FOR{$i=1,2,\dots,n$}
    \STATE{$\widehat{x}_i^{(k+1)} = (1-\omega)x_i^{(k)} + \frac{\omega}{a_{ii}} \big( b_i - \sum_{j<i} a_{ij} x_j^{(k+1)} - \sum_{j>i} a_{ij} x_j^{(k)}\big)$}
    \STATE{$x_i^{(k+1)} = \max ( \widehat{x}_i^{(k+1)} , 0 )$}
    \ENDFOR
    \ENDFOR
  \end{algorithmic}
\end{algorithm}

The objective function value sequence is nonincreasing as the iteration increases:
$V(\vecx^{(k+1)}) \leq V(\vecx^{(k)})$ (we shall give a new interpretation of this property in the next section).
This poperty is used to prove the convergence of the algorithm.

\begin{theorem}[e.g.~\cite{cr71a}]
  \label{thm:psor_spd}
  Assume that $A$ is symmetric positive definite.
  Let $\vecx^\ast$ be the unique solution of~\eqref{nqp}.
  The sequence of iterates generated by the PSOR (Algorithm~\ref{algo:sor}) converges to $\vecx^\ast$, i.e. $\vecx^{(k)} \to \vecx^\ast$ as $k\to \infty$
  if and only if $\omega\in (0,2)$.
\end{theorem}

\begin{remark}
  One may consider the following variant: 
  for $k=0,1,2,\dots$, until convergence, do
  \begin{equation}
    \begin{aligned}
       & \widetilde{\vecx}^{(k+1)} = G_{\text{SOR}} \vecx^{(k)} + \vecc_{\text{SOR}}, \\
       & \vecx^{(k+1)} = \max\{ \widetilde{\vecx}^{(k+1)}, \veczero\} ,
    \end{aligned}
    \label{naive:psor}
  \end{equation}
where $G_{\text{SOR}}$ and $\vecc_{\text{SOR}}$ are the iteration matrix and vector in the SOR method, respectively, and the $\max$ operator is taken component-wise.
Different from Algorithm~\ref{algo:sor}, the iterate here is projected 
  after all elements are updated.
  This type of projection can be incorporated with any stationary iterative method. 
  This iteration~\eqref{naive:psor} may converge to a nonnegative vector;
  however, the limit, even if it exists, depends on the relaxation parameter $\omega$ and 
  the sequence may converge to a limit which is not the unique solution of $\vecx^\ast$.
  In fact, the unique solution $\vecx^\ast$ is not necessarily a stationary point of the map, and the dissipation property $V\big(\vecx^{(k+1)}\big) \leq V\big(\vecx^{(k)}\big)$ is no longer guaranteed.
  Here, we provide two examples.
  For the problem \eqref{nqp} with
  \begin{align}
    A = 
    \begin{bmatrix}
      2 & -1 \\ -1 & 2
    \end{bmatrix},
    \quad 
    \vecb = 
    \begin{bmatrix}
      1 \\ 1
    \end{bmatrix},
  \end{align}
  the optimal solution is $\vecx^\ast = [0,1]^\trans$, which is, of course, the stationary point of Algorithm~\ref{algo:sor}.
  However, algorithm \eqref{naive:psor} converges to another vector $\lim_{k\to\infty} \vecx^{(k)} = [0,3(2-\omega)/(4-\omega)]^\trans$ even if the initial guess is set to $\vecx^{(0)} = \vecx^\ast$;
  clearly, the dissipation property is violated.
  Further, for the following case
  \begin{align}
    A = 
    \begin{bmatrix}
      2   & -1 & 0.5 \\
      -1  & 2  & -1  \\
      0.5 & -1 & 2
    \end{bmatrix},
    \quad 
    \vecb = 
    \begin{bmatrix}
      2 \\ -2 \\ 2
    \end{bmatrix},
  \end{align}
  for which the true solution is $\vecx^\ast = [0.8,0,0.8]^\trans$,  
  the algorithm \eqref{naive:psor} may stagnate.
  For example, if $\vecx^{(0)} = \veczero$, then $\vecx ^{(1)} = [\omega,0,\frac{\omega}{2}(2 - \frac{3\omega}{2} + \frac{\omega^2}{2})]^\trans$ and $\vecx^{(2)} = \veczero$ when, for example, $\omega = 1.9$.
  
\end{remark}

The next theorem summarizes the convergence property of the PSOR method applied to symmetric positive semidefinite cases.

\begin{theorem}[\cite{lu91}]
  \label{thm:psor_spsd}
  Assume that $A$ is symmetric positive semidefinite.
  The sequence of iterates that is generated by the PSOR (Algorithm~\ref{algo:sor}) converges to an element of $X^\ast$, where $X^\ast$ is the set of optimal solutions for \eqref{nqp}, if and only if $\omega \in (0,2)$.

\end{theorem}

Note that the limit of the sequence may depend on the choice of the relaxation parameter $\omega$ and the initial guess.

\begin{remark}
  The PSOR method can be applied to the nonnegative least squares problem~\eqref{nls}
  because it is equivalent to \eqref{nqp} with $A = C ^ \trans C$ and $\vecb = C^\trans \vecd$.
  However, the number of nonzero elements of $C^\trans C$ may be much larger than that of $C$, considerably increasing the computational costs.
  Instead of the SOR method, the use of the normal SOR method is recommended, which is equivalent to the SOR but avoids the explicit computation of $C^\trans C$~\cite{sa03}.
\end{remark}

\section{Adaptive projected SOR methods}
\label{sec:apsor}

In this section, we propose adaptive projected SOR (APSOR) methods for solving \eqref{nqp}.

We start our discussion by reviewing the geometric properties of the SOR method for solving $A\vecx = \vecb$:
the SOR method has a gradient structure, and the relaxation parameter can be regarded as a step size after a change of variables~\cite{ms18}.
This interpretation makes it possible to control the relaxation parameter adaptively based on a technique known in the context of continuous optimizations.
Building on our previous report~\cite{ms20}, we develop an APSOR method for solving \eqref{nqp}.
However, preliminary numerical experiments suggest that the convergence behaviour may fluctuate because of the adaptive control of the relaxation parameter; the convergence behaviour strongly depends on the initial guess, which is particularly typical for an ill-conditioned $A$ or symmetric positive semidefinite $A$.
To address these issues, we further upgrade the APSOR method in two ways.

\subsection{Geometric interpretation of the SOR method}
\label{subsec:geom_sor}

For a continuously differentiable function $V:\bbR^m\to\bbR$, 
the continuous map 
$\dg:\bbR^m \times \bbR^m \to \bbR^m$
is called a discrete gradient if it satisfies
\begin{align}
  V(\vecx) - V(\vecy) & = \dg V(\vecx,\vecy)^\trans (\vecx - \vecy), \\
  \dg V(\vecx,\vecx)  & = \nabla V(\vecx)
\end{align}
for all $\vecx, \vecy \in \bbR^m$~\cite{go96,mq99,qt96}.
The first condition is called the discrete chain rule.
The second merely requires the consistency of the discrete gradient to the gradient.
A map $\vecx^{(k)}\mapsto \vecx^{(k+1)}$ defined by
\begin{equation}
  \label{dgscheme}
  \vecx^{(k+1)} =
  \vecx^{(k)} - h P \dg V (\vecx^{(k+1)},\vecx^{(k)})
\end{equation}
with $h>0$ and a symmetric positive definite matrix $P\in\bbR^{m\times m}$ has the dissipation property
\begin{equation} \label{dissipation}
  V(\vecx^{(k+1)}) \leq V(\vecx^{(k)}),
\end{equation}
because
\begin{align}
  V(\vecx^{(k+1)}) - V(\vecx^{(k)})
   & = \dg V(\vecx^{(k+1)}, \vecx^{(k)}) ^\trans (\vecx^{(k+1)} - \vecx^{(k)})            \\
   & = -h \dg V(\vecx^{(k+1)}, \vecx^{(k)})^\trans  P \dg V(\vecx^{(k+1)}, \vecx^{(k)})  \leq 0.
\end{align}
A map of the form \eqref{dgscheme} is often referred to as the discrete gradient scheme.
Note that the symbol $h$ is called the step size because \eqref{dgscheme} can be regarded as a numerical integrator for the gradient system of the form
\begin{equation}
  \frac{\rmd}{\rmd t} \vecx(t) = -P \nabla V(\vecx(t)).
\end{equation}

The construction of discrete gradients is not unique.
For example, the discrete gradient proposed by Itoh and Abe~\cite{ia88} is defined as follows:
\begin{equation}
  \label{eq:iadg}
  \dg V(\vecx,\vecy)
  =
  \begin{bmatrix} \displaystyle
    \frac{f(x_1,y_2,\dots,y_n) - f(y_1,\dots,y_n)}{x_1-y_1}         \\
    \displaystyle
    \frac{f(x_1,x_2,y_3,\dots,y_n) - f(x_1,y_2,\dots,y_n)}{x_2-y_2} \\
    \displaystyle
    \vdots                                                          \\
    \displaystyle
    \frac{f(x_1,\dots,x_n) - f(x_1,\dots,x_{n-1},y_n)}{x_n-y_n}
  \end{bmatrix},
\end{equation}
where the case $x_i=y_i$ is interpreted as $\partial V(\vecx) / \partial x_i$.
For other definitions, see, for example,~\cite{go96,hl83,mq99,qm08}.

Let us consider the special case 
\begin{equation*}
  V(\vecx) = \frac{1}{2} \vecx^\trans A \vecx - \vecx^\trans \vecb,
\end{equation*}
where $A$ is a symmetric matrix with positive diagonal entries.
Let
$P = D^{-1}$, where $D = \diag (a_{11},a_{22},\dots,a_{nn})$.
Note that every diagonal element of $D$ is positive when $A$ is symmetric positive definite.
In this case, the Itoh--Abe discrete gradient scheme reads
\begin{align}
  x_i^{(k+1)} & = x_i^{(k)} - \frac{h}{a_{ii}} \bigg( \sum_{j<i} a_{ij} x_j^{(k+1)} + a_{ii} \frac{x_i^{(k+1)} + x_i^{(k)}}{2} + \sum_{j>i} a_{ij} x_j^{(k)} - b_i \bigg). \label{eq:iadg1} 
\end{align} 
By solving this with respect to $x_i^{(k+1)}$, we have
\begin{align}  
  x_i^{(k+1)}& = \frac{-1}{a_{ii} \big( 1+ \frac{h}{2}\big)} \bigg[ h\sum_{j<i} a_{ij}x_j^{(k+1)} - \bigg(1-\frac{h}{2}\bigg) a_{ii} x_i^{(k)} + h\sum_{j>i} a_{ij} x_j^{(k)} - h b_i\bigg]. \label{eq:iadg2}
\end{align}
The scheme~\eqref{eq:iadg2} is equivalent to the SOR method~\eqref{sor}~\cite{ms18}.

\begin{theorem}[\cite{ms18}]
  The SOR method \eqref{sor} is equivalent to the Itoh--Abe discrete gradient scheme \eqref{eq:iadg2},
  in the sense that both produce the same iterates for the same initial guess,
  when $h = 2\omega / (2-\omega)$.
\end{theorem}

Especially in the communities studying \eqref{nqp}, that the SOR method for a symmetric positive definite linear system guarantees $V(\vecx^{(k+1)})\leq V(\vecx^{(k)})$ if and only if $\omega \in (0,2)$ is a known result.
We note that 
the convergence condition in terms of $h$ is $h>0$, and
$h>0$ is equivalent to $\omega \in (0,2)$ under the relation $h = 2\omega/(2-\omega)$.
This interpretation not only suggests an alternative proof for the dissipation property $V(\vecx^{(k+1)})\leq V(\vecx^{(k)})$ but also implies a possibility of selecting the relaxation parameter $\omega$ by controlling the step size $h$.
In~\cite{ms20}, several adaptive SOR (ASOR) methods are proposed, in which the step size is adaptively controlled using ideas related to the Armijo or Wolfe conditions.

For the SOR method applied to $A\vecx = \vecb$, a direct calculation indicates that the dissipation property follows for each component-wise update.
In fact, if we focus on the update for the $i$th component, and write $\hat{x}_i^{(k+1)}$ for the update instead of $x_i^{(k+1)}$, the $i$th component of the scheme with $P=D^{-1}$ reads
\begin{align*}
  &\frac{\hat{x}_i^{(k+1)} - x_i^{(k)}}{h} \\
  & \quad = - \frac{1}{a_{ii}} \frac{V(x_1^{(k+1)}, \dots, \hat{x}_i^{(k+1)}, x_{i+1}^{(k)}, \dots, x_n^{(k)}) - V(x_1^{(k+1)}, \dots, x_{i-1}^{(k+1)}, x_i^{(k)}, \dots, x_n^{(k)})}{\hat{x}_i^{(k+1)} - x_i^{(k)}},
\end{align*} 
which simplifies to
  \begin{align}
    & V(x_1^{(k+1)}, \dots, \hat{x}_i^{(k+1)}, x_{i+1}^{(k)}, \dots, x_n^{(k)}) 
    -
  V(x_1^{(k+1)}, \dots, x_{i-1}^{(k+1)}, x_i^{(k)}, \dots, x_n^{(k)})\\
   &\quad =
  -\frac{a_{ii}}{h} (\hat{x}_i^{(k+1)} - x_i^{(k)})^2 \leq 0. 
  \label{psor:cw_dissipation}
  \end{align}
Let us consider the PSOR applied to \eqref{nqp}. 
In this case, even if the projection occurs for the update of $x_i$, a similar dissipation property
\begin{align}
  & V(x_1^{(k+1)}, \dots, x_i^{(k+1)}, x_{i+1}^{(k)}, \dots, x_n^{(k)}) 
  -
V(x_1^{(k+1)}, \dots, x_{i-1}^{(k+1)}, x_i^{(k)}, \dots, x_n^{(k)}) \\
& \quad \leq
-\frac{a_{ii}}{h} (x_i^{(k+1)} - x_i^{(k)})^2 \leq 0
\label{psor:cw_dissipation1}
\end{align}
still follows because the update can be interpreted as applying the SOR method with a nonnegative step size smaller than $h$
(note that $x_i^{(k+1)}$ is on the line segment connecting $x_i^{(k)}$ and $\hat{x}_i^{(k+1)}$).
Note that, because each component-wise update update can be understood as the SOR method with an \emph{effective} step size $h_i$ that is smaller than $h$,
\eqref{psor:cw_dissipation1} can be expressed more tightly as 
\begin{align}
  & V(x_1^{(k+1)}, \dots, x_i^{(k+1)}, x_{i+1}^{(k)}, \dots, x_n^{(k)}) 
  -
V(x_1^{(k+1)}, \dots, x_{i-1}^{(k+1)}, x_i^{(k)}, \dots, x_n^{(k)}) \\
& \quad \leq
-\frac{a_{ii}}{h_i} (x_i^{(k+1)} - x_i^{(k)})^2 \leq 0.
\end{align}
Summing up the inequality \eqref{psor:cw_dissipation1}, we see the dissipation property for the PSOR method applied to \eqref{nqp}.

\begin{lemma}
  \label{lemma1}
  Let $0<\omega<2$ or equivalently $h>0$.
  Then, for the PSOR iteration, it follows that
  \begin{equation} \label{psor:dissipation1}
    V(\vecx^{(k+1)}) - V(\vecx^{(k)})
    \leq 
    -\gamma \| \vecx^{(k+1)} - \vecx^{(k)} \|^2 
    \leq 0, 
  \end{equation}
  where 
  \begin{equation}
    \gamma = \frac{1}{h} \min (a_{11},a_{22},\dots,a_{nn}).
  \end{equation}
\end{lemma}

Theorem~\ref{thm:psor_spd} can be viewed as a consequence of the dissipation property (see appendix~\ref{appendix_1}).

\begin{remark}
  In this paper, we assume that matrix $A$ has positive diagonal entries.
  If the matrix $A$ has zero diagonal entries, the PSOR method requires preconditioning.
  Using a nonsingular matrix $M\in\bbR^{n\times n}$, we rewrite \eqref{nqp} as
  \begin{equation}
    \begin{split}
      \text{minimize} & \quad V(\vecy) = \frac{1}{2} \vecy^\trans (M^\trans AM) \vecy - \vecy^\trans (M^\trans \vecb) \\
      \text{subject to} & \quad  M\vecy \geq \veczero.
    \end{split}
  \end{equation}
  The matrix $M$ is ideally chosen such that $M^\trans A M$ has positive diagonal entries and is as sparse as possible.
  The constraint $M\vecy \geq \veczero$ indicates that the feasible region is a convex cone. 
  The projection process onto the convex cone is straightforward; it also has a low cost, especially when $M$ is sparse.
\end{remark}

\subsection{Adaptive projected SOR (APSOR) methods}

\label{subsec:apsor}

The performance of the PSOR method (Algorithm~\ref{algo:sor}) strongly depends on the choice of the relaxation parameter $\omega$,
as is the case with the SOR method applied to linear systems without any constraints.
In this subsection, we develop adaptive projected SOR (APSOR) algorithms to solve \eqref{nqp} based on the geometric interpretation discussed in section~\ref{subsec:geom_sor}.

Before presenting our algorithms, we briefly review the conventional approaches to control the relaxation parameter adaptively for linear systems.
Most existing adaptive SOR methods are categorized into two classes.
The first one is to approximate the spectral radius of the iteration matrix in each iteration using the information obtained during iterations~\cite{hy81}.
The parameter computed by this approach can be close to the optimal paramter with least iterations in practical,
but the applicability strongly depends on matrix properties.
The second one is to approximately minimize the residual during iterations~\cite{ba03} (see also~\cite{me14,rr16}).
This approach applies to any symmetric positive definite matrix, but the estimated parameter does not always satisfy the convergence condition $\omega \in (0,2)$.
Further, the computation of several additional matrix-vector multiplications is required to update the parameter,
which indicates that the cost for updating the parameter may be higher than that for updating the solution vector.

Compared to these approaches,
the adaptive SOR method that controls the step size~\cite{ms20} has the advantage that 
no additional matrix-vector multiplications are required to update the step size (or equivalently, the relaxation parameter).
Furthermore, this approach can be applied to any symmetric positive definite linear systems.
We use this idea for the PSOR method applied to \eqref{nqp}.

We start our discussion by presenting a general form of an APSOR algorithm that can build in almost all adaptive SOR techniques.
The algorithm is shown in Algorithm~\ref{algo:asor1}.

\begin{algorithm}[t]
  \caption{A general form of the APSOR method for \eqref{nqp}.}
  \begin{algorithmic}[1]
    \label{algo:asor1}
    \REQUIRE{$\omega^{(0)}\in (0,2)$, $\vecx^{(0)}\geq \veczero$}
    \FOR{$k=0,1,2,\dots$ \textbf{until} $\| \vecx^{(k+1)} - \vecx^{(k)} \| \leq \epsilon $}
    \FOR{$i=1,2,\dots,n$}
    \STATE{$\widehat{x}_i^{(k+1)} = (1-\omega^{(k)})x_i^{(k)} + \frac{\omega^{(k)}}{a_{ii}} \big( b_i - \sum_{j<i} a_{ij} x_j^{(k+1)} - \sum_{j>i} a_{ij} x_j^{(k)}\big)$}
    \STATE{$x_i^{(k+1)} = \max ( \widehat{x}_i^{(k+1)} , 0 )$}
    \ENDFOR
    \STATE{update $\omega^{(k+1)}$}
    \ENDFOR
  \end{algorithmic}
\end{algorithm}

In what follows, we understand the algorithm by replacing $\omega^{(k)}$ with $h^{(k)}$ with the relation
\begin{equation}
  h^{(k)} = \frac{2\omega^{(k)}}{2-\omega^{(k)}}.
\end{equation}
We readily see the dissipation property, which corresponds to Lemma~\ref{lemma1}. 

\begin{lemma}
  \label{lemma2}
  Assume that $h^{(k)}$ of Algorithm~\ref{algo:asor1} is bounded above, i.e. $0<h^{(k)}\leq M_h<\infty$.
  Then, for the APSOR iteration applied to \eqref{nqp}, it follows that
  \begin{equation} \label{apsor:dissipation1}
    V(\vecx^{(k+1)}) - V(\vecx^{(k)})
    \leq 
    -\tilde{\gamma} \| \vecx^{(k+1)} - \vecx^{(k)} \|^2 
    \leq 0, 
  \end{equation}
  where 
  \begin{equation}
    \tilde{\gamma} = \frac{1}{M_h} \min (a_{11},a_{22},\dots,a_{nn}).
  \end{equation}
\end{lemma}

Because of this property, global convergence is guaranteed for the symmetric positive definite cases, as long as the unwanted case $\omega^{(k)}\to 0$ or $2$ ($h^{(k)}\to 0$ or $\infty$) is avoided.

\begin{theorem}
  \label{thm:convergence2}
  Assume that the relaxation parameter $\omega^{(k)}$ in Algorithm~\ref{algo:asor1} is restricted to $\omega^{(k)} \in (\varepsilon_\omega,M_\omega)$, where $\varepsilon_\omega,M_\omega$ are predetermined constants satisfying $0<\varepsilon_\omega < M_\omega < 2$.
  When $A$ is symmetric positive definite, the sequence of iterates generated by Algorithm~\ref{algo:asor1} converges to the unique optimal solution of \eqref{nqp}.
\end{theorem}

Although the relaxation parameter, or equivalently, the step size, is not constant,
the proof goes along with the same lines as that of Theorem~\ref{thm:psor_spd} (see Appendix~\ref{appendix_1}) because of the dissipation property (Lemma~\ref{lemma2}).

\begin{remark}
  \label{rem:semi}
  Our preliminary numerical experiments suggest that the iteration seems to converge even for symmetric positive semidefinite cases.
  However, it has not been theoretically guaranteed yet.

\end{remark}

We now present an APSOR algorithm building on~\cite{ms20,rl17} and show it in Algorithm~\ref{algo:asor2}.
The idea is as follows. 
We check the Armijo condition
\begin{equation} \label{armijo}
  V(\vecx^{(k+1)}) \leq V(\vecx^{(k)}) + c_1 \nabla V(\vecx^{(k)}) ^\trans (\vecx^{(k+1)} - \vecx^{(k)})
\end{equation}
and the curvature condition
\begin{equation} \label{wolfe}
  c_2 \nabla V(\vecx^{(k)}) ^\trans (\vecx^{(k+1)} - \vecx^{(k)})
  \leq 
  \nabla V(\vecx^{(k+1)}) ^\trans (\vecx^{(k+1)} - \vecx^{(k)})
\end{equation}
after calculating $\vecx^{(k+1)}$.
If the Armijo condition is not satisfied, then the step size is decreased for the next iteration by a factor of $\rho \in (0,1)$: $h^{(k+1)} = \rho h^{(k)}$.
If the two conditions \eqref{armijo} and \eqref{wolfe} are satisfied, then the step size for the next iteration is increased by $h^{(k+1)} = \lambda_1 h^{(k)}$.
However, if the curvature condition \eqref{wolfe} is not satisfied, the current step size is deemed too small and is increased for the next iteration by a larger factor of $\lambda_2 > \lambda_1 > 1$: $h^{(k+1)} = \lambda_2 h^{(k)}$.
There are a few parameters that must be predetermined in Algorithm~\ref{algo:asor2}.
Two parameters $\varepsilon_\omega$ and $M_\omega$ are used to avoid a  situation in which $\omega^{(k)}$ takes a value close to $0$ or $2$ (because if this happens, the algorithm may stagnate).
The other parameters need to be predetermined to check the Wolfe-like conditions \eqref{armijo} and \eqref{wolfe}.
Fortunately, \cite{ms20} showed that the combination $(c_1,c_2,\lambda_1,\lambda_2,\rho) = (0.89,0.95,1.15,1.4,0.85)$ is useful for many linear systems.
This combination will also be employed in this study.
Note that our preliminary numerical experiments suggest that the convergence behaviour may not be very sensitive to the choice of these parameters, and for example, a small $c_1$ such as $c_1=10^{-4}$ as often employed in the optimization area exhibits a similar performance.

We note that in contrast to the standard applications of the Wolfe conditions that are used to determine the \emph{current} step size, the above conditions are used to determine the \emph{next} step size.
This procedure is possible because the dissipation property $V(\vecx^{(k+1)}) \leq V(\vecx^{(k)})$ always follows independently of the step size (as long as it is positive, and even if it is large).
Further, checking these conditions does not require calculation of additional matrix-vector products, which indicates that the cost for updating the step size or parameter is negligible.

\begin{algorithm}[t]
  \caption{APSOR method based on Wolfe conditions for \eqref{nqp}.}
  \begin{algorithmic}[1]
    \label{algo:asor2}
    \REQUIRE{$\vecx^{(0)}\geq \veczero$}
    \STATE{Set parameters $0<\varepsilon_\omega < M_\omega<2$, $c_1\in (0,1)$, $c_2\in (c_1,1)$, $\lambda_1>1$, $\lambda_2>\lambda_1$, $\rho\in (0,1)$}
    \STATE{$h^{(0)} = 2$}
    \STATE{$\omega^{(0)}=1$}
    \FOR{$k=0,1,2,\dots$ \textbf{until} $\| \vecx^{(k+1)} - \vecx^{(k)} \| \leq \epsilon $}
    \FOR{$i=1,2,\dots,n$}
    \STATE{$\widetilde{x}_i^{(k+1)} = \frac{1}{a_{ii}} \big( b_i - \sum_{j=1}^{i-1} a_{ij} x_j^{(k+1)} - \sum_{j=i+1}^n a_{ij} x_j^{(k)}\big)$}
    \STATE{$\widehat{x}_i^{(k+1)} = (1-\omega^{(k)}) x_i^{(k)} + \omega^{(k)} \widetilde{x}_i^{(k+1)}$}
    \STATE{$x_i^{(k+1)} = \max ( \widehat{x}_i^{(k+1)} , 0 )$}
    \IF{$V(\vecx^{(k+1)}) \leq V(\vecx^{(k)}) + c_1 \nabla V (\vecx^{(k)}) ^\trans (\vecx^{(k+1)} - \vecx^{(k)})$}
    \IF{$c_2 \nabla V (\vecx^{(k)}) ^\trans (\vecx^{(k+1)} - \vecx^{(k)}) \leq \nabla V (\vecx^{(k+1)}) ^\trans (\vecx^{(k+1)} - \vecx^{(k)})$}
    \STATE $h^{(k+1)} := \lambda_1 h^{(k)}$;
    \ELSE
    \STATE $h^{(k+1)} := \lambda_2 h^{(k)}$;
    \ENDIF
    \ELSE
    \STATE $h^{(k+1)} := \rho h^{(k)}$
    \ENDIF
    \ENDFOR
    \STATE $\displaystyle \omega^{(k+1)} := \frac{2h^{(k+1)}}{2+h^{(k+1)}}$;
    \IF{$\omega^{(k+1)} \notin (\varepsilon_\omega , M_\omega)$}
    \STATE  $h^{(k+1)} = 2$;
    \STATE  $\omega^{(k+1)} = 1$;
    \ENDIF
    \ENDFOR
  \end{algorithmic}
\end{algorithm}

\subsection{Extensions of APSOR method}
\label{subsec:var_apsor}

Preliminary experiments show the following behaviour for Algorithm~\ref{algo:asor2}.
\begin{itemize}
  \item[(a)] The convergence for Algorithm~\ref{algo:sor} looks linear, but with a slight fluctuation, whereas the rate of convergence may change during iterations.
  These are also observed for Algorithm~\ref{algo:asor2}, but the fluctuation tends to be much more significant.
  \item[(b)] The convergence performance depends on the initial guess, which is particularly typical for an ill-conditioned $A$ or symmetric positive semidefinite $A$.
\end{itemize}

Below, we present two variants of our APSOR method (Algorithm~\ref{algo:asor2}) to address these issues.

First, we consider (a).
We define the decrement for the PSOR algorithms as follows:
\begin{equation*}
  d^{(k)} := \log_{10} \| \vecx^{(k)} - \vecx^{(k-1)}\|.
\end{equation*}
The linearity of the algorithms means that $s^{(k)} := d^{(k)} - d^{(k-1)}$ is almost constant.
However, a significant deviation often occurs for Algorithm~\ref{algo:asor2}.
One method of avoiding such a situation is to fix the relaxation parameter (or equivalently, the step size) after performingx some iterations.
Below, we provide one strategy.

Let us consider the mean of $s^{(k)}$ for $m+1$ iterations, namely,
\begin{equation*}
  \overline{s}_m^{(k)} 
  = \frac{1}{m} \big(s^{(k)} + s^{(k-1)} + \cdots + s^{(k-m+1)}\big) 
  =\frac{1}{m}\big( d^{(k)} - d^{(k-m)}\big),
\end{equation*}
and the corresponding mean of the step sizes 
\begin{equation*}
  \overline{h}_m^{(k)} = \frac{1}{m+1} \big(h^{(k-1)} + h^{(k-2)} + \cdots + h^{(k-m-1)}\big)
\end{equation*}
or 
\begin{equation*}
  \overline{\omega}_m^{(k)} = \frac{1}{m+1} \big(\omega^{(k-1)} + \omega^{(k-2)} + \cdots + \omega^{(k-m-1)}\big).
\end{equation*}
Our idea is that 
if 
\begin{equation}
  \label{cond:meand}
  \overline{s}_m^{(k^\prime)} > \overline{s}_m^{(k^\prime-1)}
\end{equation}
for some $k^\prime$, 
we set the step size to
\begin{equation*}
  h^{(k)} = \overline{h}_m^{(k^\prime)}
  \quad 
  \text{or equivalently}
  \quad
  \omega^{(k)} = \overline{\omega}_m^{(k^\prime)}
\end{equation*}
for $k \geq k^\prime$.
We note that as is the case with the SOR method applied to linear systems, we need to wait for a few iterations until the convergence performance becomes almost linear.
Therefore, it seems appropriate to start checking the condition \eqref{cond:meand} after the first $l$ iterations are performed.
We recommend setting $m\approx 10$ and
determining $l$ during the iterations.
An example is $l = l^\prime + m$, where $l^\prime$ is the smallest integer such that $\overline{d}^{(l^\prime)} < -2$.  
The algorithm is summarized in Algorithm~\ref{algo:asor3}.

\begin{algorithm}[t]
  \caption{APSOR method for \eqref{nqp} that fixes the relaxation parameter during the iteration}
  \begin{algorithmic}[1]
    \label{algo:asor3}
    \REQUIRE{$\omega^{(0)}\in (0,2)$, $\vecx^{(0)}\geq \veczero$}
    \STATE{$d^{(-1)} = 0$}
    \STATE{$k=0$}
    \WHILE{$d^{(k-1)} > -2$}
    \FOR{$i=1,2,\dots,n$}
    \STATE{$\widehat{x}_i^{(k+1)} = (1-\omega^{(k)})x_i^{(k)} + \frac{\omega^{(k)}}{a_{ii}} \big( b_i - \sum_{j<i} a_{ij} x_j^{(k+1)} - \sum_{j>i} a_{ij} x_j^{(k)}\big)$}
    \STATE{$x_i^{(k+1)} = \max ( \widehat{x}_i^{(k+1)} , 0 )$}
    \ENDFOR
    \STATE{update $\omega^{(k+1)}$}
    \STATE{$k = k + 1$}
    \ENDWHILE
    \FOR{$p = 1,\dots,m+1$}
    \STATE{run lines 4--9}
    \ENDFOR
    \WHILE{$\overline{s}_m^{(k-1)} > \overline{s}_m^{(k-2)}$}
    \STATE{run lines 4--9}
    \ENDWHILE
    \STATE{$\omega = \omega^{(k+1)}$}
    \WHILE{$\| \vecx^{(k+1)} - \vecx^{(k)} \| > \epsilon $}
    \FOR{$i=1,2,\dots,n$}
    \STATE{$\widehat{x}_i^{(k+1)} = (1-\omega)x_i^{(k)} + \frac{\omega}{a_{ii}} \big( b_i - \sum_{j<i} a_{ij} x_j^{(k+1)} - \sum_{j>i} a_{ij} x_j^{(k)}\big)$}
    \STATE{$x_i^{(k+1)} = \max ( \widehat{x}_i^{(k+1)} , 0 )$}
    \ENDFOR
    \ENDWHILE
  \end{algorithmic}
\end{algorithm}

Next, we consider (b).
Finding a reasonable initial guess for the PSOR algorithms applied to NQP \eqref{nqp} is important.
Our strategy, summarized in Algorithm~\ref{algo:asor4}, involves solving the shifted problem, i.e. NQP, with the objective function 
\begin{equation*}
  V_\sigma (\vecx) = \frac{1}{2} \vecx^\trans (A+\sigma I) \vecx - \vecx^\trans \vecb,
\end{equation*}
where $\sigma >0$, using Algorithm~\ref{algo:asor2}, and set the limit (meeting convergence criteria) of the sequence to the initial guess for Algorithm~\ref{algo:asor2} applied to the original NQP problem. 
In this approach, $A+\sigma I$ is symmetric positive definite, and its condition number could be smaller than $A$.
Further, the computational cost for solving the shifted problem is unchanged from that for the original problem.
The value of $\sigma$ could substantially influence the convergence for the second line in Algorithm~\ref{algo:asor4}.
One heuristic way is to set $\sigma$ close to the minimum value of the diagonal elements of $A$.
We call this strategy pseudo-regularization.

\begin{algorithm}[t]
  \caption{APSOR method for \eqref{nqp} with the initial guess obtained by solving the shifted problem}
  \begin{algorithmic}[1]
    \label{algo:asor4}
    \STATE{run Algorithm~\ref{algo:asor2} for the shifted problem with $V_\sigma (\vecx) = \frac{1}{2} \vecx^\trans (A+\sigma I) \vecx - \vecx^\trans \vecb$ to set the initial guess $\vecx_0$}
    \STATE{run Algorithm~\ref{algo:asor2} for \eqref{nqp} using the initial guess $\vecx_0$}
  \end{algorithmic}
\end{algorithm}

\section{Numerical experiments}
\label{sec:numer}

We conducted several numerical experiments to demonstrate the efficiency of our algorithms.
First, the algorithms are tested on toy problems with a non-singular matrix $A$, where the matrix and vector $b$ are randomly set.
Second, similar experiments are conducted for toy problems with singular matrix $A$.
Finally, the algorithms are applied to image deblurring problems.  

In the numerical experiments, we use $M_\omega = 1.99$ and $\varepsilon_\omega = 0.01$ for Algorithms~\ref{algo:asor2} and~\ref{algo:asor3}.

Unless otherwise noted, we use Julia v1.6.
Codes for generating toy problems in sections~\ref{subsec:num_nonsingular} and~\ref{subsec:num_singular} are shown in Appendix~\ref{appendix_2}.
 

\subsection{Toy problems with non-singular matrix $A$}

\label{subsec:num_nonsingular}

We consider a class of sparse full rank matrices.
The size of the matrices is set to $n=10000$.
We generate four sparse symmetric positive definite matrices $A_i\in\bbR^{10000\times 10000}$ using the MATLAB function $\mathsf{sprandsym}$ with the ratio of nonzero element density $=0.1\%$.
The condition numbers of these matrices are specified as $\kappa (A_i) = 10^{3i-2}$ ($i=1,2,3,4$).
In the following numerical experiments, the termination tolerance on $\vecx$ is set as $10^{-10}$.

To observe the error behaviour of each algorithm, we need the vector $\vecb$ such that the exact solution can be attained. 
For this purpose, we first randomly prepare the exact solution $\hat{\vecx}$, almost half elements of which are set to $0$. 
Then, using a random nonnegative vector $\vecy$ satisfying $\vecx^\trans \vecy=0$, we set $\vecb = A \hat{\vecx}-\vecy$.
Results for $\kappa = 10,10^4,10^7,10^{10}$ are displayed in Figs.~\ref{fig:toy_cond1}--\ref{fig:toy_cond4}, respectively, where $\|\vecx^{(x+1)}-\vecx^{(k)}\|$, the relative error $\|\vecx^{(k)} - \hat{\vecx}\|/\|\hat{\vecx}\|$ and the relaxation parameter $\omega^{(k)}$ selected in each iteration are plotted. 
For comparison, results of the PSOR method with four choices of relaxation parameters are also plotted.

For the well-conditioned cases $\kappa=10$ (Fig.~\ref{fig:toy_cond1}) and $\kappa=10^4$ (Fig.~\ref{fig:toy_cond2}), we apply PSOR to the problems changing the relaxation parameter by increments of $0.1$, and $\omega = 1.2$ (in Fig.~\ref{fig:toy_cond1}) and $1.8$ (in Fig.~\ref{fig:toy_cond2}) are found to be the best, respectively.
The case $\kappa=10$ (Fig.~\ref{fig:toy_cond1}) is too easy, and the number of iterations required by the proposed algorithms is more than twice that required by PSOR with $\omega=1.2$ (almost optimal), but they are much faster than $\omega = 1.0$ (the Gauss--Seidel).
For the case $\kappa=10^4$ (Fig.~\ref{fig:toy_cond2}), the proposed algorithms perform well.
In particular, the convergence of APSOR that fixes the relaxation parameter after some iterations is very similar to that with $\omega=1.8$.
In fact, after 83 iterations, the relaxation parameter is fixed as $\omega = 1.77...$, which is almost optimal. 
We also note that for all cases, the performance with $\omega \in (0,1)$ was worse than that with $\omega = 1$.

For cases $\kappa=10^7$ (Fig.~\ref{fig:toy_cond3}) and $\kappa=10^{10}$ (Fig.~\ref{fig:toy_cond4}), we apply PSOR to the problem changing the relaxation parameter by increments of $0.05$, and $\omega = 1.9$ is found to be the best for both cases.
Surprisingly, APSOR performs much better than $\omega = 1.9$
(except Algorithm~\ref{algo:asor3} for the $\kappa=10^{10}$ case).

Note that changing the random seed or initial guess leads to different results. 
Although the performance is problem-dependent, the following tendencies are observed:
Algorithm~\ref{algo:asor3} shows more smooth convergence behaviour than Algorithm~\ref{algo:asor2}, although which algorithm shows faster convergence is problem-dependent;
both algorithms (Algorithms~\ref{algo:asor2} and~\ref{algo:asor3}) perform comparably to (or sometimes superior to) PSOR (Algorithm~\ref{algo:sor}) with a nearly optimal relaxation parameter.

\input{toy_cond1}
\input{toy_cond2}
\input{toy_cond3}
\input{toy_cond4}

\subsection{Toy problems with singular matrix $A$}

\label{subsec:num_singular}

We again consider a class of sparse matrices; however, this time, the matrix $A$ is symmetric positive semidefinite.
Below, we present two examples.

In the first example, we consider a matrix $A\in\bbR^{10000\times 10000}$ with $\rank(A) = 9995$.
The results of the PSOR method with several $\omega$ and Algorithm~\ref{algo:asor2} are displayed in Fig.~\ref{fig:toy_singular1}.
We observe a trend similar to that observed in section~\ref{subsec:num_nonsingular}, and APSOR performs well.

The next example considers $A\in\bbR^{100\times 100}$ with $\rank(A)=99$.
The results are displayed in Fig.~\ref{fig:toy_singular2}.
The convergence behaviour (left-top figure) shows a near stagnation for PSOR with $\omega = 1.0$ and Algorithm~\ref{algo:asor2} (note that PSOR with $\omega = 1.0$ converges after 11,834 iterations). 
A slightly similar trend is observed for PSOR with $\omega = 1.85$.
For this example, Algorithm~\ref{algo:asor4} (with pseudo-regularization) is also tested; the results are shown in the bottom figure.
Here, the parameter for the shifted problem is set to $\sigma = \min(\diag(A)) = 3030.3...$.
We observe a substantial reduction in the number of iterations:
the algorithm required only 109 iterations in total, including the number of iterations for finding the initial guess.

We applied Algorithm~\ref{algo:asor4} to other toy problems and observed a trend that the strategy with $\sigma = \min(\diag(A))$ substantially reduces the number of iterations unless each element of $A$ is too large.
Although this trend is not guaranteed theoretically (and, of course, there are counterexamples),
we stress that it is worth trying this technique when APSOR shows a sign of near stagnation behaviour.

\input{toy_singular1}

\input{toy_singular2}

\subsection{Image deblurring}

\label{subsec:id}

The above examples are artificial and are not related to any practical problems.
As a practical application, we consider image deblurring here.

Our problem was formulated as a nonnegative constrained least squares problem \eqref{nls} with the filtering operator (matrix version) $C$ and blurred and noisy image (vector version) $\vecd$.
\footnote{Problems of this type are generally formulated as the constrained least squares problem with a suitable regularization, such as the Tikhonov regularization.
This study does not employ such a regularization; thus, we can focus on the property of the original filtering operator. However, we note that the PSOR algorithms can be applied to the problem with regularization as long as the problem is of the form \eqref{nls}.}
Furthermore, we impose the upper bound for the solution: $\veczero \leq \vecx \leq \vecone$ (the PSOR algorithms can be applied to this setting in a straightforward manner).

\begin{figure}[t]
  \centering
  \includegraphics[scale=0.2]{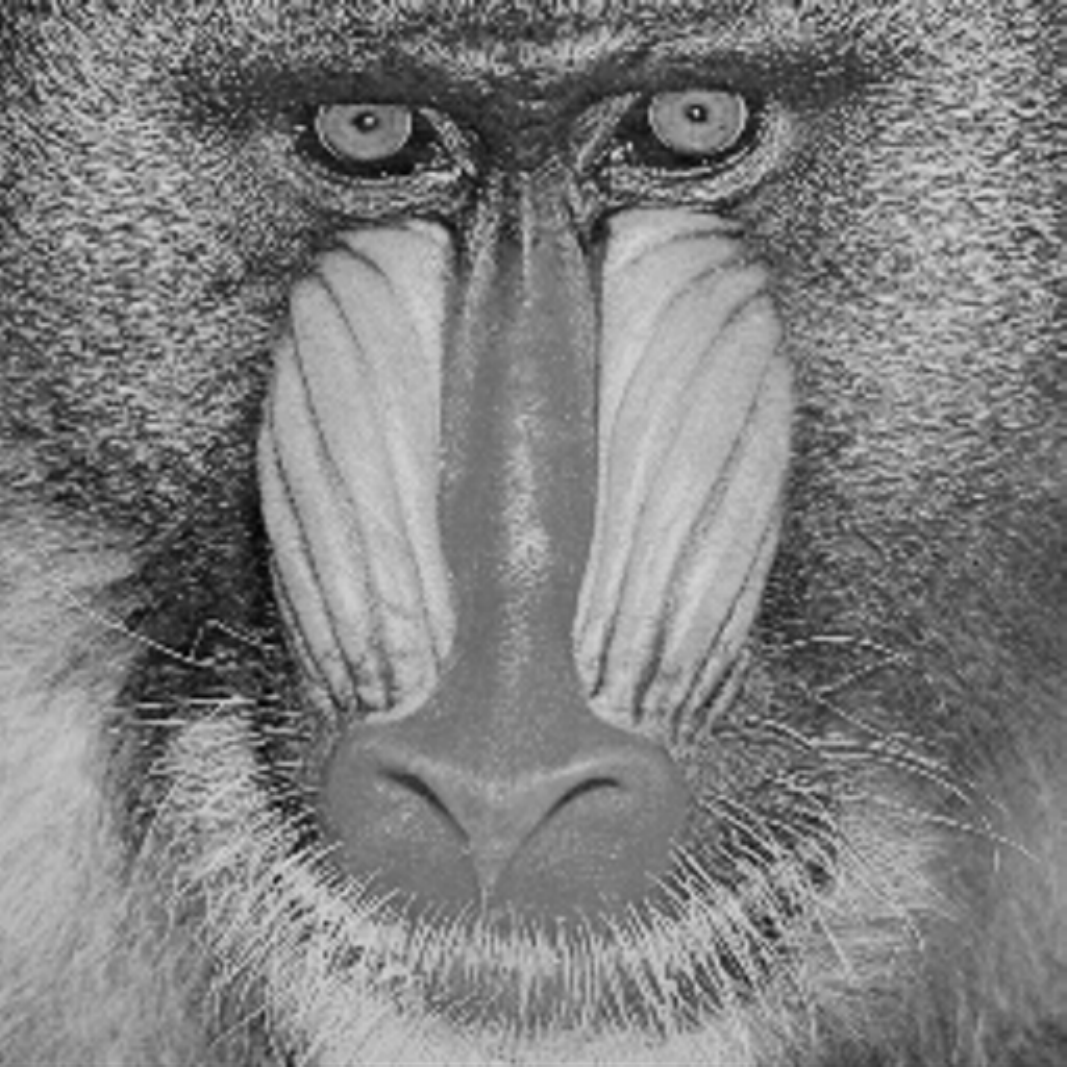}
  \quad
  \includegraphics[scale=0.2]{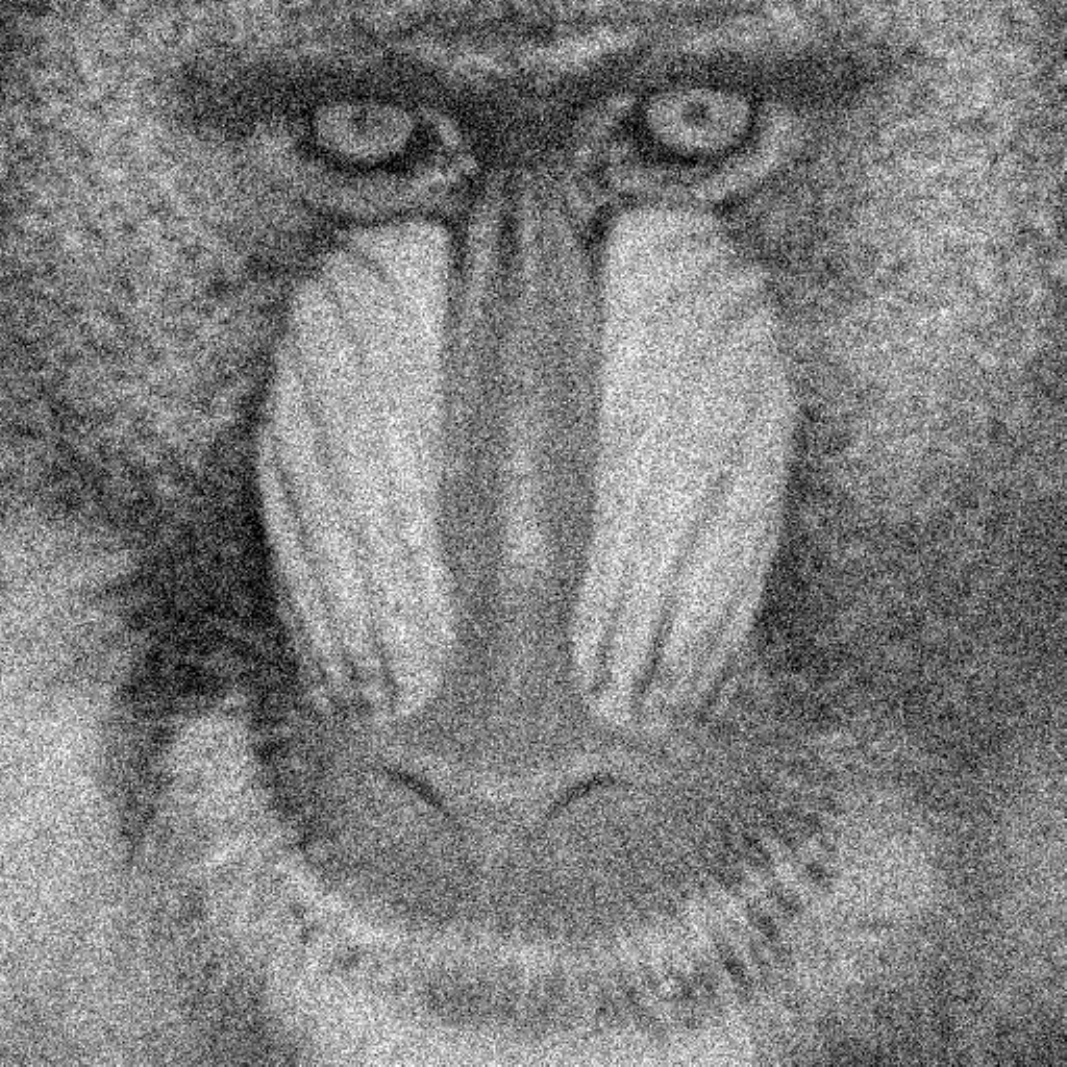}
  \quad
  \includegraphics[scale=0.2]{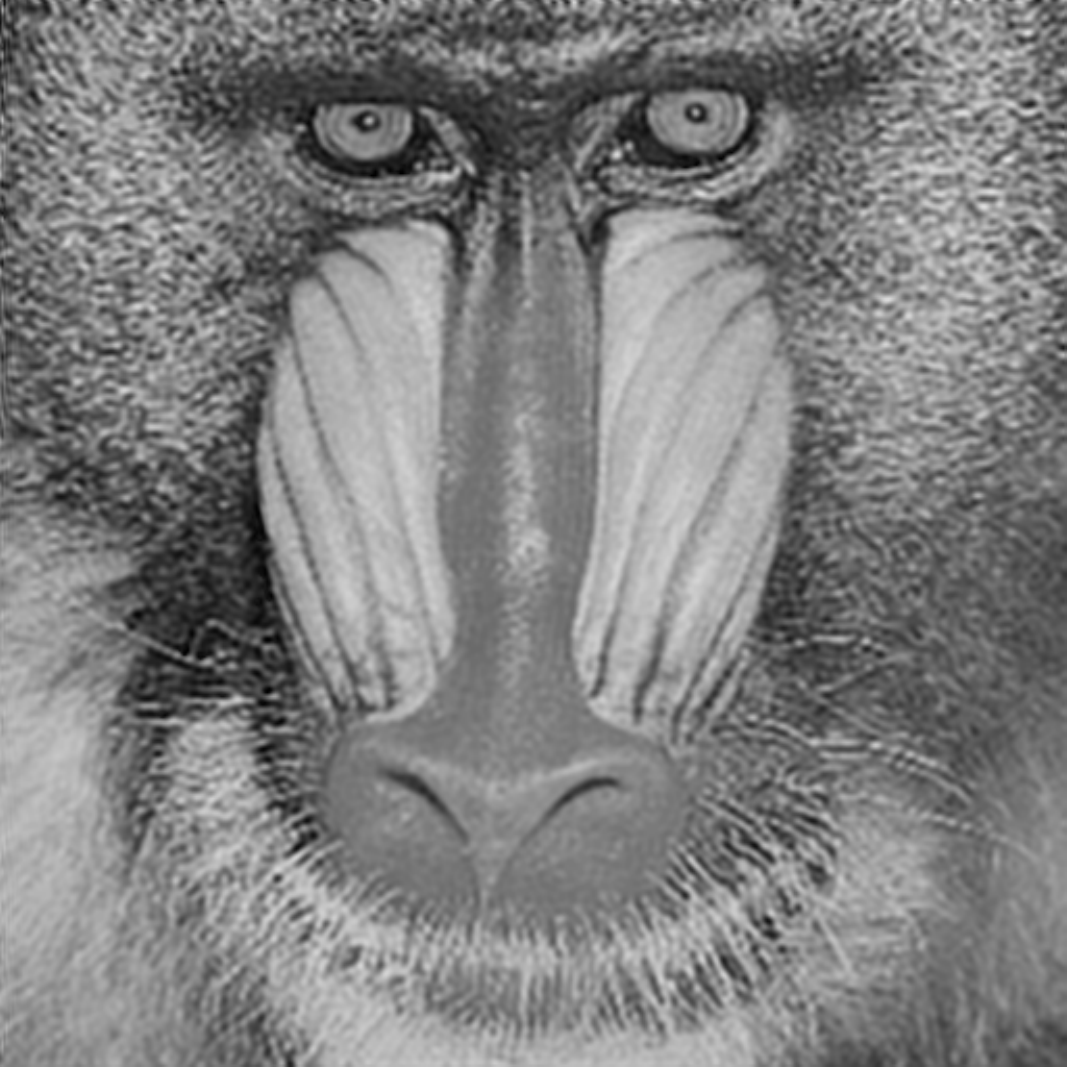}
  \caption{(Left) true baboon picture, (Center) blurred and noisy image, (Right) Deblurred image by Algorithm~\ref{algo:asor2} (50 iterations).}
  \label{fig:baboon1}
\end{figure}

In our experiment, we use a grayscale version of the mandrill, a baboon picture of the USC-SIPI image database.\footnote{The USC-SIPI Image Database, available at: \url{https://sipi.usc.edu/database/}. 
The grayscale picture is also available at Julia \textsf{TestImages} package.}
The blurred and noisy baboon picture $\vecd$ was generated by applying the Gaussian filter with length $2$ to the true image and contaminating the Gaussian noise with a standard deviation $0.1$, and is displayed in the middle picture of Fig.~\ref{fig:baboon1}
(the image size is $256\times 256$, and hence, $n=65536$).
The boundary condition is `replicate': the border pixels extend beyond the image boundaries.
Symbol $C$ denotes the filtering operator (matrix version) with this boundary condition.
Note that this matrix is symmetric positive definite; however, such a matrix can be singular depending on the filtering, size of the image, and boundary condition.
Similar experiments can be conducted even for positive semidefinite cases.

We compare Algorithm~\ref{algo:asor2} with PSOR with several choices of the relaxation parameter $\omega$.
The initial guess for all methods is set to the blurred and noisy image.
The results are shown in Fig.~\ref{fig:image1}, in which the relative error, relative residual and the relaxation parameter selected for Algorithm~\ref{algo:asor2} are displayed.
For the PSOR method, we tried $\omega = 0.1,0.2,\dots,1.9$ and show some results including the optimal case.
We observe that the best relaxation parameters in terms of relative error and the relative residual are not the same; however, note that Algorithm~\ref{algo:asor2} looks almost optimal, whereas the relaxation parameter selected for Algorithm~\ref{algo:asor2} does not remain around $0.2$ or $0.3$.

We conducted other numerical experiments with other images, and often observed the following: the superiority of Algorithm~\ref{algo:asor2} tends to be remarkable when other initial guesses are employed such as the completely black or white image;
Algorithm~\ref{algo:asor4} can be beneficial when the convergence is deemed very slow.
Note that the above experiment indicates that the proposed adaptive algorithms perform preferably even if the optimal relaxation parameter seems to be less than $1$.

\input{image1}

\section{Concluding remarks}
\label{sec:conclude}

In this study, we proposed an adaptive projected SOR method and its variants for nonnegative quadratic programming.
The relaxation parameter is controlled based on the Wolfe conditions.
The algorithms apply any nonnegative quadratic programming problems, and the cost for updating the parameter is negligible for the whole iteration.
Surprisingly, although the parameters of the Wolfe conditions need to be predetermined, the choice which was tuned for solving linear systems~\cite{ms20} works well: numerical experiments suggest that the proposed algorithms with the choice of parameters often perform comparably to (or sometimes superior to) the projected SOR method with a nearly optimal relaxation parameter, even if the optimal relaxation parameter seems to be less than $1$.

In future work, we will perform a more detailed convergence analysis for the proposed algorithms. 
In particular, as stated in Remark~\ref{rem:semi}, convergence for symmetric positive semidefinite cases needs to be clarified.


\section*{Declarations}


\begin{itemize}
\item Funding\\
This work was supported in part by JSPS, Japan KAKENHI Grant Numbers 
20H00581,20K20397,20H01822, 21K18301,
JST ACT-I Grant No. JPMJPR18US and
JST PRESTO Grant No.JPMJPR2129.
\item Conflict of interest\\
The authors declares no competing interests.
\end{itemize}

\bibliographystyle{spmpsci}      
\bibliography{references}   


\appendix
\section{Proof of Theorem~\ref{thm:psor_spd}}
\label{appendix_1}

The proof of Theorem~\ref{thm:convergence2} goes along the same lines as that of Theorem~\ref{thm:psor_spd}.
Although an equivalent statement to Theorem~\ref{thm:psor_spd} is proved in~\cite{cr71a}, we provide the proof here for the readers' convenience, incorporating some modifications and adopting the notation used in this paper.

The following lemma characterizes the unique solution to the KKT condition.

\begin{lemma}[e.g. \cite{ma77}]
  \label{lemma:ma77}
  Let $A\in\bbR^{n\times n}$.
  Let $E\in\bbR^{n\times n}$ be a positive diagonal matrix.
  Then, the following two statements are equivalent.
  \begin{enumerate}
    \renewcommand{\labelenumi}{(\roman{enumi})}
    \item $A\vecx - \vecb \geq \veczero$, \quad
          $\vecx \geq \veczero$, \quad
          $\vecx^\trans (A\vecx - \vecb) = \veczero$.
    \item $\big(\vecx - \alpha E (A\vecx - \vecb)\big)_+ - \vecx = \veczero$ holds for some or all $\alpha > 0$.
  \end{enumerate}
\end{lemma}

Theorem~\ref{thm:psor_spd} can be proved as follows.

Because $V$ is bounded from below, say by $C$, and due to the dissipation property for $V$, we see that
  \begin{equation}
    C \leq V\big(\vecx^{(k+1)}\big) \leq V\big(\vecx^{(k)}\big) \leq \cdots \leq V\big(\vecx^{(0)}\big),
  \end{equation}
  and thus, the following limit exists
  \begin{equation}
    \lim_{k\to\infty} V\big(\vecx^{(k)}\big) = V^\ast.
  \end{equation}
  Because the set 
  \begin{equation}
    V_t = \{\vecx \in \bbR^m \mid V(\vecx)\leq t\}
  \end{equation}
  is empty or compact,
  $\vecx^{(k)}$ belongs to the compact set $V_{V(\vecx^{(0)})}$.
  Therefore, for the sequence $\{\vecx^{(k)}\}_{k=0}^\infty$ there exists at least one accumulation point.

  Let $\overline{\vecx}$ be an accumulation point of $\{\vecx^{(k)}\}_{k=0}^\infty$.
  The optimal solution $\vecx^\ast$ to \eqref{nqp} satisfies the KKT condition
  \begin{equation}
    A\vecx^\ast \geq \veczero, \quad
    \vecx^\ast \geq \veczero, \quad
    (\vecx^\ast) ^\trans (A\vecx^\ast - \vecb) = \veczero.
  \end{equation}
  Note that a vector satisfying the KKT condition is unique because the objective function and feasible region of NQP are both convex.
  Below, we show that 
  \begin{equation}
    \big(\overline{\vecx} - \alpha E (A\overline{\vecx} - \vecb)\big)_+ - \overline{\vecx} = \veczero
  \end{equation}
  for some $\alpha>0$ and positive diagonal matrix $E$, which indicates due to Lemma~\ref{lemma:ma77} that
  $\overline{\vecx}$ must coincide with $\vecx^\ast$.

  Let $\{\vecx^{(k_l)}\}_{l=0}^\infty$ be a subsequence that converges to $\overline{\vecx}$.
  Because of the continuity of $V$, we see that
  \begin{equation}
    \lim_{l\to\infty} V\big(\vecx^{(k_l)}\big) = V(\overline{\vecx}).
  \end{equation}
  From the convergence of $\{V(\vecx^{(k)})\}_{k=0}^\infty$ and \eqref{psor:dissipation1}, it follows that 
  \begin{equation}
    0 = \lim_{l\to\infty} \Big( V\big(\vecx^{(k_l+1)}\big) - V\big(\vecx^{(k_l)}\big)\Big)
    \leq - \lim_{l\to\infty} \gamma \big\| \vecx^{(k_l+1)} - \vecx^{(k_l)} \big\|^2 \leq 0.
  \end{equation}
  This relation and the compactness of $V_{V(\vecx^{(0)})}$ lead to
  \begin{align}
    0 & = \lim_{l\to\infty} \big\| \vecx^{(k_l+1)} - \vecx^{(k_l)}\big\|^2
    = \lim_{l\to\infty} \big\| \big(\vecx^{(k_l+1)} - \overline{\vecx}\big) - \big(\vecx^{(k_l)} - \overline{\vecx}\big) \big\|^2 \\
      & = \lim_{l\to\infty} \big\| \vecx^{(k_l+1)} - \overline{\vecx} \big\|^2 .
  \end{align}
  Therefore,
  \begin{equation}
    \lim_{l\to\infty} \vecx^{(k_l+1)} = \lim_{l\to\infty} \vecx^{(k_l)} = \overline{\vecx}.
  \end{equation}

  It follows that
  \begin{align}
    0 & = \lim_{l\to\infty} \big\| \vecx^{(k_l+1)} - \vecx^{(k_l)} \big\|^2 
    = \lim_{l\to\infty} \sum_{i=1}^{m} \big(x_i^{(k_l+1)} - x_i^{(k_l)}\big)^2 \\
      & = \lim_{l\to\infty} \sum_{i=1}^{m}
    \bigg\{ \bigg[
    \frac{-1}{a_{ii} \big( 1+ \frac{h}{2}\big)}
    \bigg( h\sum_{j<i} a_{ij}x_j^{(k_l+1)} - \Big(1-\frac{h}{2}\Big) a_{ii} x_i^{(k_l)} + h\sum_{j>i} a_{ij} x_j^{(k_l)} - h b_i\bigg) 
    \bigg]_{+} - x_i^{(k_l)} \bigg\}^2                                         \\
      & =  \sum_{i=1}^{m}
    \bigg\{ \bigg[
      \frac{-1}{a_{ii} \big( 1+ \frac{h}{2}\big)}
      \bigg( h\sum_{j<i} a_{ij}\overline{x}_j - \Big(1-\frac{h}{2}\Big) a_{ii} \overline{x}_i + h\sum_{j>i} a_{ij} \overline{x}_j - h b_i\bigg)
    \bigg]_{+} - \overline{x}_i \bigg\}^2                                      \\
      & = \sum_{i=1}^{m}
    \bigg\{ \Big[
      \overline{x}_i - \frac{h}{1+\frac{h}{2}} \frac{1}{a_{ii}} (A\overline{\vecx} - \vecb)_i 
      \Big]_{+} - \overline{x}_i \bigg\}^2.
  \end{align}
  Therefore, by setting
  \begin{equation}
    \alpha = \frac{h}{1+\frac{h}{2}},\quad
    E = D^{-1},
  \end{equation}
  where $D$ is the diagonal matrix of $A$,
  we get from Lemma~\ref{lemma:ma77} that $\overline{\vecx} = \vecx^\ast$.
  This shows that any convergent subsequence converges to the unique solution $\vecx^\ast$; hence, the whole sequence also converges to $\vecx^\ast$.
  This completes the proof.

  \section{Codes for generating toy problems}
  \label{appendix_2}

  The vector $\vecb$ for all toy problems are generated by the following Julia code:

  \begin{verbatim}
      Random.seed!(1)
      xtrue = max.(randn(n),0)
      idx = findall(xtrue.==0)
      ytrue = zeros(n)
      ytrue[idx] = abs.(randn(length(idx)))
      dot(xtrue,ytrue)
      b = A*xtrue - ytrue
  \end{verbatim}

Four matrices used in Section~\ref{subsec:num_nonsingular} are generated by the following MATLAB codes:
  \begin{verbatim}
      rng(1)
      n = 10000;
      density = 0.001; 
      cond_num = 1e1; % (or 1e4, 1e7,1e10) 
      A = sprandsym(n,density,linspace(1,cond_num,n));
  \end{verbatim}

  The matrix of the first example in Section~\ref{subsec:num_singular} is generated by the following MATLAB code:

  \begin{verbatim}
      rng(1)
      n = 10000;
      density = 0.005;
      cond_num = 1e10;
      r = 9995;
      a = zeros(1,n-r-1);
      b = linspace(0,cond_num,r+1);
      a = [a,b];
      A = sprandsym(n,density,a);
  \end{verbatim}

  and the second is generated by

  \begin{verbatim}
      rng(1)
      n = 100;
      density = 0.1;
      cond_num = 1e5;
      A = sprandsym(n,density,linspace(0,cond_num,n));
  \end{verbatim}

\end{document}